\documentclass[12pt,oneside,reqno]{amsart}
\usepackage{txfonts}
\usepackage{bbm}

\usepackage{amsmath}
\usepackage{graphicx}
\usepackage{mathrsfs}
\usepackage{stmaryrd}
\usepackage{amsfonts}

\usepackage{enumerate,amsmath,amssymb,amsthm,booktabs}

\numberwithin{equation}{section}

\newcommand{\be}{\begin{eqnarray}}
\newcommand{\ee}{\end{eqnarray}}
\newcommand{\ce}{\begin{eqnarray*}}
\newcommand{\de}{\end{eqnarray*}}
\newtheorem{theorem}{Theorem}[section]
\newtheorem{lemma}[theorem]{Lemma}
\newtheorem{remark}[theorem]{Remark}
\newtheorem{definition}[theorem]{Definition}
\newtheorem{proposition}[theorem]{Proposition}
\newtheorem{Examples}[theorem]{Example}
\newtheorem{corollary}[theorem]{Corollary}

\def\eps{\varepsilon}

\def\p{\partial}

\def\[{{\Big[}}
\def\]{{\Big]}}
\def\<{{\langle}}
\def\>{{\rangle}}
\def\({{\Big(}}
\def\){{\Big)}}

\def\bx{{\mathbf{x}}}

\def\dif{{\mathord{{\rm d}}}}

\def\no{\nonumber}
\def\={&\!\!=\!\!&}
\def\bt{\begin{theorem}}
\def\et{\end{theorem}}
\def\bl{\begin{lemma}}
\def\el{\end{lemma}}
\def\br{\begin{remark}}
\def\er{\end{remark}}
\def\bd{\begin{definition}}
\def\ed{\end{definition}}
\def\bp{\begin{proposition}}
\def\ep{\end{proposition}}
\def\bc{\begin{corollary}}
\def\ec{\end{corollary}}
\def\bx{\begin{Examples}}
\def\ex{\end{Examples}}

\def\mE{{\mathbb E}}

\def\mN{{\mathbb N}}

\def\mP{{\mathbb P}}

\def\mR{{\mathbb R}}

\def\bP{{\mathbf P}}
\def\bQ{{\mathbf Q}}

\def\1{{\mathbf{1}}}

\def\sF{{\mathscr F}}

\def\sL{{\mathscr L}}

\def\sR{{\mathscr R}}

\def\geq{\geqslant}
\def\leq{\leqslant}

\def\bP{{\mathbf P}}

\allowdisplaybreaks

\begin{document}

\title{H\"older estimates for nonlocal-diffusion equations with drifts}

\date{}

\author{{Zhen-Qing Chen} \quad
and \quad {Xicheng Zhang}}

\address{Zhen-Qing Chen:
Department of Mathematics, University of Washington, Seattle, WA 98195, USA\\
Email: zqchen@uw.edu
 }
\address{Xicheng Zhang:
School of Mathematics and Statistics, Wuhan University,
Wuhan, Hubei 430072, P.R.China\\
Email: XichengZhang@gmail.com
 }

\thanks{The research of ZC is partially supported
by NSF grant DMS-1206276.
The research of XZ is partially supported by NNSFC grant of China (Nos. 11271294, 11325105).}

\begin{abstract}
We study a class of nonlocal-diffusion equations with drifts, and derive a priori $\Phi$-H\"older estimate for the solutions by using a purely probabilistic argument, where $\Phi$ is an intrinsic scaling function for the equation.

\bigskip



\end{abstract}

\maketitle

\section{Introduction}

It is well-known that a priori H\"older estimates play a key role in the study of non-linear equations.
There is a wealth literature on this for differential operators.
In this paper, we study a priori H\"older estimates of parabolic functions for a large class of time-dependent non-local operators with (time-dependent) gradient drifts  on $\mR^d$
$$
\sL^b_t u =\sL_t u+b_t\cdot\nabla u.
$$
Here $b_t(x)$ is an $\mR^d$-valued measurable function on $[0, \infty) \times\mR^d$ and
$\sL_t$ is a time-dependent purely non-local operator given by
\begin{equation}\label{e:1.1}
\sL_t u(x)=\int_{\mR^d}[u(x+z)+u(x-z)-2u(x)]\kappa_t(x,z)\dif z,
\end{equation}
where $\kappa_t(x,z)$ is symmetric in $z$ (i.e. $\kappa_t(x,z)=\kappa_t(x,-z)$)
 and satisfies
\begin{align}\label{Con1}
\sup_{t,x}\int_{\mR^d}(1\wedge|z|^2)\kappa_t(x,z)\dif z\leq C_1,
\end{align}
and for some regularly varying function $\phi$ with index $\alpha\in[0,2]$ (see Definition \ref{Def31} below),
\begin{align} \label{Kappa0}
\frac{c_1}{\phi(|z|)|z|^d}\leq\kappa_t(x,z)\leq& \frac{c_1^{-1}}{\phi(|z|)|z|^d}
\quad \hbox{for }  |z|\leq 3.
\end{align}
Note that under condition \eqref{Con1}, $\sL_t u(x)$ is well-defined
and bounded for every $u\in C^2_b(\mR^d)$.
Since $\kappa_t (x, z)$ is symmetric in $z$, we can rewrite $\sL_t u(x)$
in \eqref{e:1.1} as
\begin{eqnarray*}
\sL_t u(x)&=& 2 {\rm p.v.} \int_{\mR^d} \left( u(x+z)-u(x)\right) \kappa_t (x, z)
\dif z \\
&=& 2 \int_{\mR^d}  \left( u(x+z)-u(x)-\nabla u(x) \cdot z \1_{|z|\leq 1}\right) \kappa_t (x, z) \dif z
\end{eqnarray*}
for $u\in C^2_b (\mR^d)$.

When $\kappa_t (x, z)= c |z|^{-d-\alpha}$ for some suitable constant $c>0$
and $\alpha \in (0, 2)$,
$\sL_t$ is just the usual fractional Laplacian $\Delta^{\alpha/2}$.
Recently  Silvestre \cite{Si1} proved the following a priori H\"older estimate for the fractional-diffusion equation with drift
\begin{align}
\p_t u=\Delta^{ {\alpha}/{2}}u+b_t\cdot\nabla u . \label{Fr}
\end{align}
There are constants $C>0$ and $\beta\in(0,1)$ such that for any classical solution $u$ of equation (\ref{Fr}),  $x, y\in \mR^d$ with $|x-y|\leq 1$ and $0<s\leq t\leq 1$,
\begin{align}
|u(t,x)-u(s,y)|\leq
 C \| u \|_{L^\infty ([0, 1]\times \mR^d)} \,
\frac{|x-y|^\beta+|t-s|^{\beta/\alpha}}{t^{\beta/\alpha}},  \label{Hol}
\end{align}
provided $b$ in the H\"older class $C^{1-\alpha}$ when $\alpha\in(0,1)$, and  bounded measurable when $\alpha\in[1,2)$,
where the constants $C$ and $\beta$ only depend on $d,\alpha$, $\|b\|_\infty$, as well as on the $(1-\alpha )$-H\"older norm
 of $b$ when $\alpha\in(0,1)$.
See \cite{EP} for recent result on local regularity of solutions to
$$
\Delta^{ {\alpha}/{2}}u+b (x) \cdot\nabla u =f
$$
in Sobolev spaces with $\alpha \in (0, 1)$.

In the literature, if $\alpha\in(1,2)$, the equation (\ref{Fr}) is usually referred to as the {\it subcritical} case since the fractional Laplacian
$\Delta^{{\alpha}/{2}}$ is of higher order than that of the gradient part $b_t \cdot \nabla$.
For $\alpha=1$, it is called the {\it critical} case since the fractional Laplacian has the same order as the first order gradient term. For $\alpha\in(0,1)$,
it is known as the {\it supercritical} case because the fractional Laplacian is of lower order than the drift term, and the drift term can be stronger than the diffusion term in small scales.
This explains why one needs $b$ to be H\"older continuous for the above a priori estimate in the supercritical case. It should be noted that the following scaling property plays a crucial role in paper \cite{Si1}: for $\lambda>0$, let
$u^\lambda(t,x):=\lambda^{-\alpha}u(\lambda^\alpha t, \lambda x)$
and $b^\lambda(t,x):=b(\lambda^\alpha t, \lambda x)$,
then $u^\lambda$ satisfies
$$
\p_t u^\lambda=\Delta^{{\alpha}/{2}}u^\lambda+\lambda^{\alpha-1}b^\lambda\cdot\nabla u^\lambda.
$$
We mention that the H\"older's estimate (\ref{Hol}) has been used in the study of well-posedness of multidimensional critical Burger's equations in
Zhang \cite{Zh1}.

In this paper, we are concerned with the following nonlocal-diffusion equation with drift $b$:
\begin{align}\label{Eq1}
\p_t u =\sL^b_t u  =\sL_t u+b_t\cdot\nabla u.
\end{align}
Following \cite{Ka-Mi}, we define
\begin{align}\label{Phi0}
\Phi(r):=\left(\int^2_r\frac{\dif s}{s \phi(s)}\right)^{-1},\ \ r\in(0,1),
\end{align}
which is a continuous increasing function.
The purpose of introducing this function $\Phi$ is to deal with the case
when $\phi$ in \eqref{Kappa0} is a regularly varying function of order $0$.
It is known  (see \eqref{Phi2} below) that
$\lim_{r\to 0}\frac{\Phi(r)}{\phi(r)}=\alpha$.
Thus $\phi$ and $\Phi$ are comparable when $\alpha \in (0, 2]$
so we could use $\phi$ in place of $\Phi$ in this case.
The function $\Phi$ will be used to measure the modulus of continuity
of parabolic functions; see Theorem \ref{Main}.
Such a result would be trivial if $\lim_{r\to 0} \Phi ( r ) >0$.
Thus without loss of generality,
 we will assume in this paper that
 \begin{equation} \label{AS1}
 \lim_{r\to 0} \Phi ( r ) =0.
 \end{equation}
 Note that
$$
\lim_{r\to 0}\Phi(r)=0\iff \int^2_0\frac{\dif s}{s\phi(s) }=\infty\stackrel{(\ref{Kappa0})}{\iff} \int_{\mR^d}\kappa_t(x,z)\dif z=\infty.
$$
In \cite{Ka-Mi}, Kassmann and Mimica called $\Phi$ an intrinsic scaling function
and obtained   $\Phi$-H\"older regularity of harmonic functions
for the (time-independent) non-local operator $\sL_0$
with $b_0=0$ under conditions (\ref{Con1}), (\ref{Kappa0}) and \eqref{AS1}.
When $\kappa_0 (x, z)$ satisfies \eqref{Kappa0} for all $z\in \mR^d$
and $b_0=0$,  H\"older regularity of harmonic functions of $\sL_0$ was first established in Bass and Levin \cite{Ba-Le}.

In this paper, we use ``:=" as a way of definition.
Let $\mR^+:=[0, \infty)$, and $C^{\infty}_c
(\mR^+\times \mR^d)$ the space of smooth functions
 with compact support in $\mR^+\times \mR^d $.
A probability measure $\bQ$ on the Skorokhod space
${\mathbb D}([0, \infty); \mR^+\times \mR^d)$ is said to
 be a solution to the martingale problem for
$(\sL_t^b, C^{\infty}_c(\mR^+\times \mR^d))$
with initial value $(t, x)\in \mR^+ \times \mR^d$
if   $\bQ(Z_0=(t, x))=1$ and for every
$f\in C^{\infty}_c(\mR^+\times \mR^d)$,
\begin{equation}\label{e:1.7}
 M^f_s:=f(s+t, X_s) -f(t, X_0)-\int_0^s (\partial_r+\sL_{t+r}^b) f(t+r, X_r) \dif r
\end{equation}
is a $\bQ$-martingale.
The martingale problem for $(\sL^b_t, C^\infty_c(\mR^+\times \mR^d))$
with initial value $(t, x)\in \mR^d$
is said to be well-posed if it has a unique solution.

Throughout this paper, we assume conditions \eqref{Con1} and \eqref{Kappa0} and
\eqref{AS1} hold. Here are our main results.

\bt\label{Main}
Assume that  the martingale problem for $(\sL^b_t, C^\infty_c (\mR^+\times \mR^d))$ is well-posed for every initial value
$(t, x)\in \mR^+\times \mR^d$.  Let $\Phi$ be defined by (\ref{Phi0}).

\begin{enumerate}[\rm (i)]
\item Suppose that \,
$\displaystyle \liminf_{r\to 0}r/\Phi(r)=0$
and for some $C_2>0$,
\begin{align}
\frac{\Phi(r)}{\Phi(s)}\leq C_2\frac{r}{s}  \quad \hbox{for }
  \ 0<s\leq r\leq 1, \label{EE01}
\end{align}
and $b_t (x)$ is continuous in $x$ for each $t>0$ having
$$
\|b/(1+|x|)\|_\infty:=\sup_{t\in [0, 1], x\in \mR^d}|b_t(x)|/(1+|x|)<\infty
$$
with
\begin{align}
|b_t(x)-b_t(y)|\leq C_b|x-y|/\Phi(|x-y|)  \quad \hbox{for }\  |x-y|\leq 1 \label{bb0}
\end{align}
for some $C_b>0$.
Then there are constants $\beta\in(0,1)$ and $\lambda=\lambda (\|b/(1+|x|)\|_\infty)>0$ such that for any classical solution $u$ of (\ref{Eq1}),
and for any $0<t\leq t_0\leq 1$ and $|x_0-x|+ \lambda (t_0-t) \leq\Phi^{-1}(t_0)$,
\begin{align}
|u(t_0,x_0)-u(t,x)|\leq 16
\|u\|_{L^\infty ([0, 1]\times \mR^d)}\,
 t_0^{-\beta}((t_0-t)+\Phi(|x_0-x|+\lambda(t_0-t)))^\beta.\label{Es0}
\end{align}

\item
Suppose \,
$\displaystyle \liminf_{r\to 0}r/\Phi(r)>0$
and $b$ is a bounded measurable function on $\mR^+\times \mR^d$.
Then there is a constant $\beta\in(0,1)$ such that for any classical solution $u$ of (\ref{Eq1}),
and for any $0<t\leq t_0\leq 1$ and $|x_0-x|\leq\Phi^{-1}(t_0)$,
\begin{align}
|u(t_0,x_0)-u(t,x)|\leq 16
\|u\|_{L^\infty ([0, 1]\times \mR^d)}\,
t_0^{-\beta}((t_0-t)+\Phi(|x_0-x|))^\beta.\label{Es00}
\end{align}
\end{enumerate}
\et
\br
\rm
Condition (\ref{EE01}) is automatically satisfied by (\ref{Pr}) below when $\phi$ (also $\Phi$) is regularly varying with $\alpha\in[0,1)$.
Moreover, condition \eqref{EE01} is satisfied when $\Phi$ is comparable to a convex function $\Psi$ with $\Psi (0)=0$
as in this case $r/\Psi (r)$ is increasing.
The following table gives some examples of regularly varying functions
$\phi$ that satisfy conditions \eqref{Con1}, \eqref{Kappa0} and \eqref{AS1}, where the $\alpha$ stands for the index of regularly varying function $\phi$, the last column
denotes the modulus of continuity in spatial variable $x$
required for $\mR^d$-valued function $b_t (x)$.

\medskip

$$
\begin{tabular}{c|cccc}\toprule
Case &$\alpha$& $\phi(s)$&$\Phi(r)$&$b$\\ \midrule
(i) & $0$& $\ln \frac{3}{s}$ & $\asymp(\ln\ln \frac{3}{r})^{-1}$ & $s\ln\ln \frac{3}{s}$\\ \midrule
(i) & $0$& $1$ & $(\ln \frac{2}{r})^{-1}$ & $s\ln \frac{2}{s}$\\ \midrule
(i) & $0$& $(\ln \frac{3}{s})^{-1}$ & $\asymp(\ln \frac{3}{r})^{-2}$ & $s(\ln \frac{3}{s})^2$\\ \midrule
(i) &$(0,1)$& $s^\alpha$ & $\asymp r^{\alpha}$&$s^{1-\alpha}$\\ \midrule
(i) &\ \ \ $1$\ \ \  &$\ \ \ s(\ln \frac{3}{s})^\beta,\ \beta>0\ \ \ $ &\ \ \  $\asymp r(\ln \frac{3}{r})^\beta$\ \ \ &$\ \ \ (\ln\frac{3}{s})^{-\beta}$\ \ \ \\ \midrule
(ii) &$1$ &$s(\ln \frac{3}{s})^\beta,\ \beta\leq 0$ & $\asymp r(\ln \frac{3}{r})^\beta$&$1$\\ \midrule
(ii) &$(1,2)$& $s^\alpha$ & $\asymp r^{\alpha}$&$1$\\ \midrule
(ii) &$2$ &$s^2(\ln\frac{3}{s})^{\beta},\ \beta>1$ & $\asymp r^2(\ln\frac{3}{r})^{\beta}$&$1$\\
\bottomrule
\end{tabular}
$$
 \er

\medskip

\br
\rm When $\phi (r)=r^\alpha$ with $\alpha \in (0, 2)$ (and so $\Phi (r)$ is comparable to $r^\alpha$), condition
$\displaystyle \liminf_{r\to 0}r/\Phi(r)=0$
corresponds precisely to the supercritical case
$0<\alpha<1$, while $\displaystyle \liminf_{r\to 0}r/\Phi(r)>0$ corresponds to
$1\leq \alpha <2$. So Theorem \ref{Main} contains the main results of
Silvestre \cite{Si1}
as a particular case.
Unlike \cite{Si1}, in the supercritical case (i), we do not need to assume  $b$ is bounded. Moreover,
in this paper we can not only deal with
more general  but also time-dependent nonlocal operator $\sL^b_t$.
In particular, taking $b=0$, we obtain  a priori $\Phi$-H\"older estimates
for parabolic functions of time-dependent non-local operator $\sL_t$.
It contains as a special case   a priori $\Phi$-H\"older estimates
for harmonic functions $u(x)$ of  non-local operator $\sL_0$, which
is the main result of
Kassmann and Mimica \cite{Ka-Mi} when $b_0=0$. \qed
\er

\br \rm
In this paper, we concentrate on a priori H\"older estimates for parabolic functions.
For results on the well-posedness of the martingale problems for $\sL^b_t$,
we refer the reader to \cite{Ab-Ka, Ch-Wa, Ch-Zh} and the references therein.
We remark here that if $\kappa_t (x, z)$ is independent of $(t, x)$,
$\sL_t=\sL_0$ is the generator of a L\'evy process $Y$. When $b_t (x)$ is uniformly
Lipschitz in $x$, it is easy to show that for every initial data $(t, x)\in \mR^+\times \mR^d$, stochastic differential equation
\begin{equation}\label{e:1.13}
 \dif X^t_s = \dif Y_s + b_{t+s} (X^t_s) \dif s
\end{equation}
has a unique strong solution $X^t$ with $X^t_0=x$.
(This can be done as follows. For each $\omega \in \Omega$,
ODE $dZ^t_s = b_{t+s} (Z^t_s +Y_s(\omega)) \dif s$ with $Z^t_0 (\omega) =x$
has a unique solution. Then $X^t_s := Z^t_s+Y_t$ is the unique solution to
SDE \eqref{e:1.13} with $X^t_0=x$.)
Hence in this case, the martingale problem for $(\sL, C^\infty_c (\mR^+\times \mR^d))$ is well-posed and so Theorem \ref{Main} is applicable.
It is important to note that the $\Phi$-H\"older estimate in Theorem \ref{Main} does not depend on the Lipschitz constant of $b_t (x)$.

When the martingale problem for $(\sL^b_t, C^\infty_c (\mR^+\times \mR^d))$ is well-posed for every initial value $(t, x)\in \mR^+\times \mR^d$, there is a
space-time Hunt process $Z_s=(V_0+s, X_s)$ having $\partial_s + \sL^v_{V_0+s}$
as its infinitesimal generator. For any bounded classical solution $u$ of (\ref{Eq1}), by It\^o's formula,  $u(s+t, X_s)$ is a $\bP_{(t, x)}$-martingale for every
$(t, x)\in \mR^+\times \mR^d$. This is the only property of $u$ we used in the proof of Theorem \ref{Main}. Hence the conclusion of Theorem \ref{Main}
holds for any bounded function $u$ on $\mR^+\times \mR^d$ such that $u(s+t, X_s)$ is a $\bP_{(t, x)}$-martingale for every $(t, x)\in \mR^+\times \mR^d$. \qed
\er

The approach of this paper is purely probabilistic.
Our tool is the time-inhomogeneous strong Markov process $X$
 determined by the solution
of the martingale problem for $(\sL^b_t, C^\infty_c (\mR^+\times \mR^d))$.
In Section 2, we prove an abstract result on H\"older's continuity in terms of certain estimates on exiting and hitting probabilities, which is motivated by the approaches in \cite{Ba-Le, Ch-Ku}.
This probabilistic approach has its origin in Krylov and Safanov \cite{KS}
for diffusion processes associated with second order non-divergence form
differential operators.
In Section 3, we prove our main results by verifying the abstract conditions. The L\'evy system of the strong Markov process $X$
plays a key role in establishing these exiting and hitting probabilities.

\section{An abstract criterion for H\"older's regularity}

Let $\Omega$ be the space of   c\`adl\`ag functions from $\mR^+ =[0,\infty)$ to $\mR^d$, which is  endowed with the Skorokhod topology.
Let $X_s(\omega)=\omega_s$ be the coordinate process over $\Omega$.
Define the space-time process
$$
Z_s:=(V_s, X_s),\ \ V_s:=V_0+s.
$$
Let $\{\sF_s^0; s\geq 0\}$ be the  natural filtration generated by $X$. Suppose $\{\mP_{(t,x)}; t\geq 0, x\in \mR^d\}$ is a family of probability measures over $(\Omega,\sF^0_\infty)$
 so that $Z=(\Omega,\sF^0_\infty,\sF^0_s, Z_s, \mP_{(t,x)})$ is a time-{\it homogenous} strong Markov processes with state space $\mR^+\times\mR^d$ with
$$
\mP_{(t,x)}\big(Z_0=(t,x)\big)=1.
$$
Denote by $\{\sF_s: s\geq 0\}$ the minimal augmented filtration of $Z$.
Note that under $\mP_{(t, x)}$, $\{X^t_s:=X_{s+t}; s\geq 0\}$
is a possibly {\it time-inhomogeneous} strong Markov process with
$$
\mP_{(t,x)}(X^t_s=x,\ s\in[0,t])=1.
$$

For a Borel set $A\subset\mR^+\times\mR^d$, denote by $\sigma_A,\tau_A$ the hitting time and exit time of $A$, i.e.,
$$
\sigma_A:=\inf\{s\geq 0: Z_s\in A\},\quad \tau_A:=\inf\{s\geq 0: Z_s\notin A\}.
$$
\bd
A nonnegative Borel measurable function $u(t,x)$ on $\mR^+\times\mR^d$ is called $Z$-harmonic (or simply parabolic) in a relatively open subset $D$ of $\mR^+\times\mR^d$
if for each relatively compact open subset $A\subset D$ and every $(t,x)\in A$,
\begin{align}
u(t,x)=\mE_{(t,x)}[u(Z_{\tau_A})].\label{ET0}
\end{align}
\ed
\br{\rm\label{Le14}
Condition (\ref{ET0}) is equivalent to the following. For any ($\sF_s$)-stopping time $\tau$,
\begin{align}
u(t,x)=\mE_{(t,x)}\[u(Z_{\tau\wedge\tau_A})\].\label{H8}
\end{align}
Indeed, assume that (\ref{ET0}) holds for any $(t,x)\in A$. In view of $Z_\tau\in A$ on $\{\tau<\tau_A\}$, we have
$$
\1_{\tau<\tau_A}u(Z_\tau)=\1_{\tau<\tau_A}\mE_{Z_\tau}[u(Z_{\tau_A})].
$$
Let $\{\theta_t; t\geq 0\}$ be the usual shift operators on $\Omega$. By the strong Markov property, we have
\begin{align*}
\mE_{(t,x)}\[\1_{\tau<\tau_A}u(Z_\tau)\]&=\mE_{(t,x)}\[\1_{\tau<\tau_A}\mE_{Z_{\tau}}[u(Z_{\tau_A})]\]
=\mE_{(t,x)}\[\1_{\tau<\tau_A}\mE_{(t,x)}\[u(Z_{\tau_A}\circ\theta_\tau)|\sF_\tau\]\]\\
&=\mE_{(t,x)}\[\1_{\tau<\tau_A}u(Z_{\tau_A}\circ\theta_\tau)\]\quad
\hbox{since }  \{\tau<\tau_A\}\in\sF_\tau .
\end{align*}
Since $\tau_A=\tau+\tau_A\circ\theta_\tau$ on $\{\tau<\tau_A\}$ and $Z_{\tau_A}\circ\theta_\tau=Z_{\tau+\tau_A\circ\theta_\tau}$, we obtain
$$
\mE_{(t,x)}\[\1_{\tau<\tau_A}u(Z_\tau)\]=\mE_{(t,x)}\[\1_{\tau<\tau_A}u(Z_{\tau_A})\],
$$
which implies that
\begin{align*}
\mE_{(t,x)}\[u(Z_{\tau\wedge\tau_A})\]&=\mE_{(t,x)}\[\1_{\tau<\tau_A}u(Z_\tau)\]+\mE_{(t,x)}\[\1_{\{\tau\geq\tau_A\}}u(Z_{\tau_A})\]
=\mE_{(t,x)}\[u(Z_{\tau_A})\]=u(t,x).
\end{align*}
}\er

\vspace{5mm}
Let $\Phi:(0,2)\to[0,\infty)$ be a continuous and strictly increasing function
with $\Phi (1)=1$.
 Write for $r>0$,
$$
  B(r):=\{z\in\mR^d: |z|<r\} \quad \hbox{and}  \quad
 Q(r):=[0,\Phi(r))\times B(r).
$$
Define
\begin{align}
\varphi_a(r):=\Phi^{-1}(a\Phi(r))\mbox{ for $a>0$}, \quad
D_a(r):=Q(\varphi_a(r))\setminus Q(\varphi_{\sqrt{a}}(r)) \mbox{ for $a>1$}.\label{Da}
\end{align}
Notice that
$$
\varphi_1(r)=r \quad \mbox{and} \quad
a\mapsto\varphi_a(r)\mbox{ is strictly increasing.}
$$
We make the following assumptions:
\begin{enumerate}[{\bf (H$_1$)}]
\item There exist constants $C_{3}, C_{4}\geq 1$ such that for each $a\geq C_{4}$ and $r,R\in(0,1)$ with $\varphi_a(r)\leq R$,
\begin{align}
\sup_{(t_0,x_0)\in Q(r)}\mP_{(t_0,x_0)}\Big(X_{\tau_{Q(r)}}\notin B(R)\Big)\leq C_{3}\frac{\Phi(r)}{\Phi(R)}.\label{ERY3}
\end{align}
\item There is an increasing sequence of positive numbers $\{a_k; k\geq 1\}\subset (1, \infty)$ with $\lim_{k\to \infty} a_k=\infty$ such that
    for every $a\in  \{a_k; k\geq 1\}$,
there exists a constant
$\gamma_a\in (0,a]$
so that for each $r\in(0,1)$ with $\varphi_a(r)\leq 1$, there is a radon measure $\mu_r$ over $D_a(r)$
such that for any compact subset $K\subset D_a(r)$ with $\mu_r(K)\geq\frac{1}{3}\mu_r(D_a(r))$,
\begin{align}
\inf_{(t_0,x_0)\in Q(r)}\mP_{(t_0,x_0)}\Big(\sigma_K<\tau_{Q(\varphi_a(r))}\Big)\geq \frac{\gamma_a}{a},\label{H3}
\end{align}
and $\lim_{k\to \infty} \gamma_{a_k}=\infty$.
\end{enumerate}

\br{\rm If $X_s$ is a continuous process, then {\bf (H$_1$)} is automatically satisfied.
}\er

\bt\label{Th24}
Under {\bf (H$_1$)} and {\bf (H$_2$)}, there exists a constant $\beta\in(0,1)$, which only depends on $C_{3}, C_{4}$ and $\gamma_a$, such that for each $r\in(0,1)$,
every bounded measurable function $u$ on
$[0, 1] \times\mR^d$ that is parabolic in $Q(r)$,
\begin{align}
|u(t,x)-u(0,0)|\leq 8\left(\frac{t\vee\Phi(|x|)}{\Phi(r)}\right)^\beta
\|u\|_{L^\infty ([0, \Phi (r)]\times \mR^d)} \quad \hbox{for }
(t,x)\in Q(r).\label{RE}
\end{align}
\et

\begin{proof}
Our proof is adapted from
Chen and Kumagai \cite[Theorem 4.14]{Ch-Ku}.
Fix $r\in (0, 1)$. Without loss of generality, we may assume $0\leq u\leq 1$
on $[0, \Phi (r)]\times \mR^d$.
Otherwise, instead of $u$, we may consider
$$
\tilde u_t (x) =\frac{u_t (x)-\inf_{(s, y)\in  [0, \Phi (r)]\times \mR^d}u_s (y)}{\sup_{(s, y)\in  [0, \Phi (r)]\times \mR^d}u_s (y)-\inf_{(s, y)\in  [0, \Phi (r)]\times \mR^d}u_s (y)}.
$$
\\
{\bf (i)}  Define for $n\in\mN$,
$$
r_n:=\varphi_{a^{1-n}}(r),\quad  s_n:=2 b^{1-n},
$$
where $a>1$ from $\{a_k; k\geq 1\}$ and $b\in (1, 2)$ to
 be determined below.
Observe that $\Phi (r_n)=a\Phi (r_{n+1})$.
 Clearly,
$$
\varphi_a(r_{n+1})=r_{n}\mbox{ and } r_n\downarrow 0, \ \ s_n\downarrow 0.
$$
For simplicity of notation, we write
$$
Q_n:=Q(r_n),\quad   M_n:=\sup_{Q_n}u,\quad m_n:=\inf_{Q_n}u.
$$
We are going to prove that the oscillation of $u$ over $Q_k$
\begin{align}
\hbox{osc}_{Q_k} u:=M_k-m_k\leq s_k, \quad  k\in\mN.\label{ER0}
\end{align}
If this is proven, then (\ref{RE}) follows. In fact, for any $(t,x)\in Q_1$, there is an $n\in\mN$ such that
$$
(t,x)\in Q_n\setminus Q_{n+1},
$$
which means that
$$
\Phi(r_{n+1})\leq t<\Phi(r_n)=a \Phi (r_{n+1})
\quad \hbox{or}\quad  r_{n+1}\leq|x|<r_n.
$$
In this case, we have
$$
|u(t,x)-u(0,0)|\leq M_n-m_n\leq s_n=2b a^{-n\ln b/\ln a}\leq 2b\left(\frac{t\vee\Phi(|x|)}{\Phi(r)}\right)^{\frac{\ln b}{\ln a}},
$$
and (\ref{RE}) follows with $\beta=\ln b/\ln a$.
\\
\\
{\bf (ii)} We now prove (\ref{ER0}) by an inductive argument. First of all, clearly,
$$
M_1-m_1\leq 1\leq s_1=2,\ \ M_2-m_2\leq 1\leq s_2=2/b.
$$
Next suppose that $M_k-m_k\leq s_k$ for all $k=1,\cdots, n$. Define
$$
A:=\Big\{z\in D_a(r_{n+1}): u(z)\leq \tfrac{m_n+M_n}{2}\Big\}.
$$
By considering $1-u$ instead of $u$ if necessary, we may assume that
$$
\mu_{r_{n+1}}(A)\geq\tfrac{1}{2}\mu_{r_{n+1}}(D_a(r_{n+1})),
$$
where $\mu_{r_{n+1}}$ is given in {\bf (H$_2$)}.
(Note here we are interested in the oscillation $\hbox{osc}_{Q_k} u=M_k-m_k$ not
on the exact values of $M_k$ and $m_k$.)
Since $\mu_{r_{n+1}}$ is regular, there is a compact subset $K\subset A$ such that
\begin{align}
\mu_{r_{n+1}}(K)\geq\tfrac{1}{3}\mu_{r_{n+1}}(D_a(r_{n+1})).\label{We4}
\end{align}
For any $\eps>0$, let us choose $z_1, z_2\in Q_{n+1}$ so that
$$
u(z_1)\leq m_{n+1}+\eps,\ \ u(z_2)\geq M_{n+1}-\eps.
$$
If one can show
\begin{align}
u(z_2)-u(z_1)\leq s_{n+1},\label{ER21}
\end{align}
then
$$
M_{n+1}-m_{n+1}-2\eps\leq s_{n+1}\Rightarrow M_{n+1}-m_{n+1}\leq s_{n+1},
$$
and (\ref{ER0}) is thus proven.
\\
\\
{\bf (iii)} Now, we show (\ref{ER21}). Since $z_2\in Q_{n+1}\subset Q_n$, if we define $\tau_n:=\tau_{Q_n}$, then by (\ref{H8}) we have
\begin{align}
&u(z_2)-u(z_1) =\mE_{z_2}\Big[u(Z_{\tau_{n}\wedge\sigma_K})-u(z_1)\Big] \no\\
=& \left( \mE_{z_2}\Big[u(Z_{\sigma_K})-u(z_1); \sigma_K<\tau_{n}\Big] +\mE_{z_2}\Big[u(Z_{\tau_{n}})-u(z_1); \sigma_K\geq\tau_{n};Z_{\tau_{n}}\in Q_{n-1}\Big]\right) \no\\
& +\mE_{z_2}\Big[u(Z_{\tau_{n}})-u(z_1); \sigma_K\geq\tau_{n}, Z_{\tau_{n}}\notin Q_{n-1}\Big]\no \\
=:&I_1+I_2.\label{We6}
\end{align}
For $I_1$, since $u(z_1)\geq m_{n+1}\geq m_n\geq m_{n-1}$, by the inductive hypothesis we have
\begin{align}
I_1&\leq \Big(\tfrac{m_n+M_n}{2}-m_n\Big)\mP_{z_2}(\sigma_K<\tau_{n})+(M_{n-1}-m_{n-1})\mP_{z_2}(\sigma_K\geq\tau_{n})\no\\
& \leq \tfrac{s_n}{2}\mP_{z_2}(\sigma_K<\tau_{n})+s_{n-1}(1-\mP_{z_2}(\sigma_K<\tau_{n}))\no\\
&\leq s_{n-1}(1-\mP_{z_2}(\sigma_K<\tau_{n})/2)\leq s_{n+1}b^2(1-\gamma_a/(2a)),\label{We5}
\end{align}
where the last step is due to (\ref{We4}) and ({\bf H$_2$}).
For $I_2$, we similarly have
\begin{align*}
I_2&=\sum_{i=1}^{n-2}\mE_{z_2}\Big[u(Z_{\tau_{n}})-u(z_1); \sigma_K\geq\tau_{n}, Z_{\tau_{n}}\in Q_{n-i-1}\setminus Q_{n-i}\Big]\\
&\qquad+\mE_{z_2}\Big[u(Z_{\tau_{n}})-u(z_1); \sigma_K\geq\tau_{n}, Z_{\tau_{n}}\notin Q_1\Big]\\
&\leq\sum_{i=1}^{n-2}s_{n-i-1}\mP_{z_2}\Big(Z_{\tau_{n}}\notin Q_{n-i}\Big)+\mP_{z_2}\Big(Z_{\tau_{n}}\notin Q_1\Big).
\end{align*}
Noticing that
$$
\mP_{z_2}\Big(Z_{\tau_{n}}\notin Q_{n-i}\Big)=
\mP_{z_2}\Big(X_{\tau_{Q( r_n})}\notin B ( r_{n-i})\Big),
$$
by {\bf (H$_1$)}, we further have for
$a > \max \{C_4, b\}$,
$$
I_2\leq C_{3}\sum_{i=1}^{n-2}s_{n-i-1}\frac{\Phi(r_{n})}{\Phi(r_{n-i})}+C_{3}
\frac{\Phi(r_{n})}{\Phi(r_1)}
=2C_{3}b^{2-n}\sum_{i=1}^{n-2}(b/a)^{i}+C_{3}a^{1-n}
\leq s_{n+1}b^2\left(\frac{C_{3}b}{a-b}+\frac{C_{3}}{2a}\right),
$$
which, together with (\ref{We6}) and (\ref{We5}), yields that
$$
u(z_2)-u(z_1)\leq s_{n+1} b^2\left(1-\frac{\gamma_a}{2a}+\frac{C_{3}b}{a-b}+\frac{C_{3}}{2a}\right)
\leq s_{n+1} b^2 \left(1-\frac{\gamma_a}{3a} \right) \leq s_{n+1}
$$
provided we take $a=a_k$ large enough and $b$ close to $1$ as $\lim_{k\to \infty} \gamma_{a_k}=\infty$.
This completes the proof.
\end{proof}

\section{Proof of Theorem \ref{Main}}

We first recall the definition and properties of  regularly varying functions.
\bd\label{Def31}
A measurable and positive function $\phi: (0,2)\to (0,\infty)$ is said to vary regularly at zero with index $\alpha\in\mR$ if for every $\lambda>0$,
$$
\lim_{r\to 0}\frac{\phi(\lambda r)}{\phi(r)}=\lambda^\alpha.
$$
We call such $\phi$
a regularly varying function. All regularly varying functions with index $\alpha$ is denoted by $\sR_\alpha$.
\ed

We list some properties of $\phi\in\sR_\alpha$ for later use (cf. \cite[p.25-28]{Bi-Go-Te} and \cite{Ka-Mi}).

\bp\label{Pro1}
Let $\alpha\geq 0$
and $\phi\in\sR_\alpha$  be bounded away from $0$ and $\infty$ on any compact subset of $(0,2)$.
For any $\delta>0$, there is a constant $C_{5}=C_{5}(\delta)\geq 1$
such that for all $r,s\in(0,1]$,
\begin{align}
\frac{\phi(r)}{\phi(s)}\leq C_{5}\max\left\{\Big(\frac{r}{s}\Big)^{\alpha+\delta}, \Big(\frac{r}{s}\Big)^{\alpha-\delta}\right\},\label{Pr}
\end{align}
and for any $\beta>\alpha-1$,
\begin{align}
\lim_{r\to 0}\frac{\phi(r)}{r^{\beta+1}}\int^r_0 \frac{s^{\beta}}{\phi(s)}\dif s&=(\beta-\alpha+1)^{-1},\label{Pr1}\\
\lim_{r\to 0}r^{\beta+1-\alpha}\phi(r)\int^2_r \frac{1}{s^{\beta+2-\alpha}\phi(s)}\dif s&=(\beta-\alpha+1)^{-1}.\label{Pr11}
\end{align}
Moreover, if we define
\begin{align}
\Phi(r):=\left(\int^2_r\frac{1}{\phi(s)s}\dif s\right)^{-1},\label{Phi}
\end{align}
then $\Phi\in\sR_\alpha$ and
\begin{align}\label{Phi2}
\lim_{r\to 0}\frac{\Phi(r)}{\phi(r)}=\alpha.
\end{align}
In particular, (\ref{Pr}) and (\ref{Pr1}) also hold for $\Phi$, and for some $C_{6}> 1$,
\begin{align}\label{Phi1}
\phi (2s)\leq C_{6} \phi (s) \quad \hbox{and} \quad
\Phi(2s)\leq C_{6}\Phi(s) \quad \hbox{ for } s\in(0,1/2).
\end{align}
\ep

\medskip
We now return to the setting in Section 1.
By normalizing the function $\phi$ in \eqref{Kappa0} by a constant multiple,
we may and do assume the scale function $\Phi$ defined by
\eqref{Phi0} has the property that $\Phi (1)=1$.
Consider the nonlocal operator $\sL^b_t $ in (\ref{Eq1}). We assume
\begin{itemize}
\item [{\bf (MP)}] The martingale problem for $(\sL^b_t, C^\infty_c (\mR^+\times \mR^d))$ is well-posed for every initial value $(t, x)\in \mR^+\times \mR^d$.
\end{itemize}

Denote by $\mP_{(t,x)}$ the law of the unique solution to the martingale problem for $(\sL^b_t, C^\infty_c (\mR^+\times \mR^d))$ with initial value $(t, x)\in \mR^+\times \mR^d$.
By \cite[Theorems 4.3.12 and 4.4.2]{Et-Ku}), $Z_s=(V_0+s, X_s)$ is a Hunt process with $\mP_{(t,x)}(V_0=t \hbox{ and } X_0=x)=1$ and so it has a L\'evy system
that describes the jumps of $Z$.
By a similar argument as that for \cite[Theorem 2.6]{CKS}, we have the following.

\bt\label{T:l2}
 Assume {\bf (MP)} holds. Then for any $(t,x)\in \mR^+\times\mR^d$ and any non-negative measurable function $f$ on $\mR^+
\times \mR^d\times \mR^d$ vanishing on $\{(s, x, y)\in \mR^+ \times
\mR^d\times \mR^d: x=y\}$ and $(\sF_t)$-stopping time $T$,
\begin{equation}\label{e4}
\mE_{(t,x)} \left[\sum_{s\leq T} f(s,X_{s-}, X_s) \right]= \mE_{(t,x)} \left[
\int_0^T \left( \int_{\mR^d} f(s,X_s, y) \kappa_{s+t}(X_s, y-X_s)\dif y \right) \dif s \right].
\end{equation}
\et

Next we prove the following estimate, which implies {\bf (H$_1$)}.
\bl\label{Le34}
Let $C_{6}$ be as in (\ref{Phi1}). Under (\ref{Con1}), (\ref{Kappa0}) and {\bf (MP)},
there is a constant $C_{7}\geq 1$ such that for all $a\geq C_{6}$ and $r,R\in(0,1)$ with $\varphi_a(r)\leq R$,
$$
\sup_{(t_0,x_0)\in Q(r)}\mP_{(t_0,x_0)}\Big(X_{\tau_{Q(r)}}\notin B(R)\Big)\leq C_{7}\frac{\Phi(r)}{\Phi(R)},
$$
where $\Phi$ is defined by (\ref{Phi}).
\el
\begin{proof}
For simplicity of notation, we write $z=(t_0,x_0)$.
Note that $r<\varphi_a(r)$ so we have by formula (\ref{e4}),
\begin{align*}
&\mP_{z}\Big(X_{\tau_{Q(r)}}\notin B(R)\Big)=\mE_{z}\left(\sum_{0<s\leq\tau_{Q(r)}}\1_{\{X_{s-}\in B(r), X_s\in B(R)^c\}}\right)
=\mE_{z}\int^{\tau_{Q(r)}}_0\!\!\!\!\int_{B(R)^c}\kappa_{s+t}(X_s, X_s-y)\dif y\dif s\\
&=\mE_{z}\int^{\tau_{Q(r)}}_0\!\!\!\!\int_{B(2)\cap B(R)^c}\kappa_{s+t}(X_s, X_s-y)\dif y\dif s
+\mE_{z}\int^{\tau_{Q(r)}}_0\!\!\!\!\int_{B(2)^c}\kappa_{s+t}(X_s, X_s-y)\dif y\dif s=:I_1+I_2.
\end{align*}
By (\ref{Phi1}), we have
$\varphi_a(r)\geq 2r$ for $a\geq C_{6}$,
which implies that
$$
 |x-y|\geq |y|-|x|\geq |y|/2 \quad \hbox{for } x\in B(r) \hbox{ and } y\in B(R)^c\subset B(\varphi_a(r))^c .
$$
For $I_1$, by (\ref{Kappa0}) and (\ref{Pr}) we have
\begin{align*}
I_1&\leq\mE_{z}\int^{\tau_{Q(r)}}_0\!\!\!\int_{B(2)\cap B(R)^c}\frac{c_1^{-1}}{\phi(|X_s-y|)|X_s-y|^d}\dif y\dif s\\
&\leq C\mE_{z}\tau_{Q(r)}\int_{B(2)\cap B(R)^c}\frac{\dif y}{\phi(|y|)|y|^d}\leq C\mE_{z}\tau_{Q(r)}/\Phi(R).
\end{align*}
On the other hand, by (\ref{Con1}) we clearly have
\begin{align*}
I_2\leq\mE_{z}\int^{\tau_{Q(r)}}_0\!\!\!\!\int_{B(1)^c}\kappa_{s+t}(X_s, y)\dif y\dif s\leq C_1\mE_{z}\tau_{Q(r)}.
\end{align*}
Hence, by (\ref{AS1}),
\begin{align*}
\mP_{z}\Big(X_{\tau_{Q(r)}}\notin B(R)\Big)\leq\mE_{z}\tau_{Q(r)}\Big(C_1+C/\Phi(R)\Big)\leq C_7\mE_{z}\tau_{Q(r)}/\Phi(R),
\end{align*}
which yields the desired estimate by $\tau_{Q(r)}\leq\Phi(r)$.
\end{proof}

Before verifying {\bf (H$_2$)}, we need the following lemma.
\bl\label{Le56}
Let $\Phi$ be defined by (\ref{Phi}). Suppose that one of the following conditions holds:
\begin{enumerate}[\rm (i)]
\item $\liminf_{r\to 0}r/\Phi(r)=0$ and for some $C_2>0$,
\begin{align}
\frac{\Phi(r)}{\Phi(s)}\leq C_2\frac{r}{s},\quad  0<s\leq r\leq 1,\label{EE1}
\end{align}
and
for some $C_b>0$,
\begin{align}
|b_t(x)|\leq C_b|x|/\Phi(|x|), \quad  |x|\leq1.\label{bb}
\end{align}
\item $\liminf_{r\to 0}r/\Phi(r)>0$ and $b$ is bounded measurable.
\end{enumerate}
Then there exists a constant $C_{8}\geq 1$ such that for all $r\in(0,1)$, $x_0\in B(r)$ and $t_0\in[0,1]$,
\begin{align}
\mP_{(t_0,x_0)}(\tau_{B(x_0,r)}<t)\leq \frac{C_{8} t}{\Phi(r)},\quad t>0.\label{ERY3b}
\end{align}
In particular, for any $a\geq C_{6}^2>1$ and $r\in(0,1)$ with $\varphi_a(r)\leq 1$,
\begin{align}\label{We6b}
\sup_{x_0\in B(r)}\mP_{(t_0,x_0)}
\Big(\tau_{Q(\varphi_{a}(r))}<
 \Phi(r)\Big)\leq\frac{C_{8}}{\sqrt{a}},
\end{align}
where $C_6$ is the positive constant in (\ref{Phi1}).
\el

\begin{proof}
Given $f\in C^2_b(\mR^d)$ with $f(0)=0$ and $f(x)=1$ for $|x|\geq 1$,    set
$$
f_r(x):=f((x-x_0)/r),\ \ r>0.
$$
By the optional stopping theorem,
\begin{align}
\mP_{(t_0,x_0)}\Big(\tau_{B(x_0,r)}<t\Big)\leq \mE_{(t_0,x_0)}f_r\Big(X_{\tau_{B(x_0,r)}\wedge t}\Big)
=\mE_{(t_0,x_0)}\int^{\tau_{B(x_0,r)}\wedge t}_0\sL^b_{s+t_0} f_{r}(X_{s})\dif s.\label{ERY2}
\end{align}
On the other hand, by the definition of $\sL^b_s$ and (\ref{Kappa0}), we have
\begin{eqnarray*}
|\sL^b_s f_r(x)|&=& \left|\int_{\mR^d}(f_r(x+z)+f_r(x-z)-2f_r(x))\kappa_s(x,z)\dif z+b_s(x)\cdot\nabla f_r(x)\right|\\
&\leq & C\int_{|z|\leq r}\frac{\|\nabla^2 f_r\|_\infty}{\phi(|z|)|z|^{d-2}}\dif z+C\int_{1\geq |z|\geq r}\frac{\|f_r\|_\infty}{\phi(|z|)|z|^{d}}\dif z
+\|f_r\|_\infty\int_{|z|\geq 1}\kappa_s(x,z)\dif z +\|\nabla f_r\|_\infty|b_s(x)|\\
&\leq & \frac{C}{r^2}\int^r_0\frac{s\dif s}{\phi(s)}+C\int^1_r\frac{\dif s}{\phi(s)s}+C+\frac{C|b_s(x)|}{r}\\
&\leq & \frac{C}{\phi(r)}+\frac{C}{\Phi(r)}+C+\frac{C|b_s(x)|}{r}
\qquad \hbox{by \eqref{Pr1}   and   \eqref{Pr11}.}
\end{eqnarray*}
Substituting this into (\ref{ERY2}) and using (\ref{Phi2}) and (\ref{AS1}), we obtain
\begin{align}
\mP_{(t_0,x_0)}\Big(\tau_{B(x_0,r)}<t\Big)\leq \frac{Ct}{\Phi(r)}
+\mE_{(t_0,x_0)}\int^{\tau_{B(x_0,r)}\wedge t}_0\frac{C|b_{s+t_0}(X_s)|}{r}\dif s.\label{ERY4}
\end{align}
In case (i), since
$$
\mbox{$|x|\leq |x-x_0|+|x_0|\leq 2r$ for $x\in B(x_0,r)$ and $x_0\in B(r)$,}
$$
(\ref{ERY3b}) follows by (\ref{ERY4}) and
$$
\frac{|b_{s+t_0}(x)|}{r}\stackrel{(\ref{bb})}{\leq} \frac{C_b|x|}{r\Phi(|x|)}\stackrel{(\ref{EE1})}{\leq} \frac{2C_bC_2}{\Phi(2r)}\leq \frac{C}{\Phi(r)}.
$$
In case (ii), (\ref{ERY3b}) follows by (\ref{ERY4}) and $|b_{s+t_0}(X_s)|\leq
\| b\|_{L^\infty(\mR^+\times \mR^d)}$
as well as $\frac{1}{r}\leq \frac{C}{\Phi(r)}$.
\\
\\
On the other hand, by (\ref{Phi1}), we have for any $a\geq C_{6}^2$,
$$
\Phi(2\varphi_{\sqrt{a}}(r))\leq C_{6} \sqrt{a}\Phi(r)\leq a\Phi(r),
$$
which implies that for $x_0\in B(r)$ and $x\in B(x_0,\varphi_{\sqrt{a}}(r))$,
$$
|x|\leq |x-x_0|+|x_0|\leq \varphi_{\sqrt{a}}(r)+r\leq 2\varphi_{\sqrt{a}}(r)\leq\varphi_{a}(r).
$$
Hence
$B(x_0,\varphi_{\sqrt{a}}(r))\subset B(\varphi_a (r))$
and
$$
\mP_{(t_0,x_0)}\Big(\tau_{Q(\varphi_a(r))}<\Phi(r)\Big)
\leq \mP_{(t_0,x_0)}\Big(\tau_{B(x_0,\varphi_{\sqrt{a}}(r))}
< \Phi(r)\Big)
\stackrel{(\ref{ERY3b})}{\leq}\frac{C_{8}\Phi(r)}{\Phi(\varphi_{\sqrt{a}}(r))}=\frac{C_{8}}{\sqrt{a}}.
$$
The proof is complete.
\end{proof}
For $a>1$ and $r\in(0,1)$ with $\varphi_a(r)\leq 1$, let $D_a(r)$ be defined by (\ref{Da}).
Define a measure
$$
\mu_r(A):=\int^{\Phi(r)}_0\!\!\!\int_{\mR^d}\frac{1_A(s,y)\Phi(|y|)}{\phi(|y|)|y|^d}\dif y\dif s,\ \ A\subset D_a(r).
$$
The above definition of $\mu_r$ arises   naturally when estimating the lower bound of
$\mP_{(t_0,x_0)}\Big(\sigma_K<\tau_{Q(\varphi_a(r))}\Big)$  using the L\'eve system of $Z$; see
\eqref{e:3.16} below.
Clearly, we have
\begin{eqnarray}
\mu_r(D_a(r))&=& \Phi(r)\int_{\varphi_{\sqrt{a}}(r)\leq |y|\leq\varphi_a(r) }\frac{\Phi(|y|)}{\phi(|y|)|y|^d}\dif y
=\omega_d\Phi(r)\int_{\varphi_{\sqrt{a}}(r)}^{\varphi_a(r) }\frac{\Phi(s)}{\phi(s)s}\dif s \nonumber \\
&=& \omega_d \Phi (r) \int_{\varphi_{\sqrt{a}}(r)}^{\varphi_a(r) } \frac{1}{\Phi (s)} d\Phi (s) =\frac{1}2 \omega_d\Phi(r)\ln a ,\label{We7}
\end{eqnarray}
where $\omega_d$ is the sphere area of the unit ball.

\bl\label{Le36}
 Suppose (\ref{Con1}), (\ref{Kappa0}), {\bf (MP)} and the assumptions of Lemma \ref{Le56} hold.  There exist $a_0\geq 1$ and $c_2\in(0,1)$
such that for each $a\geq a_0$ and $r\in(0,1)$ with $\varphi_a(r)\leq 1$, and any compact subset $K\subset D_a(r)$ with $\mu_r(K)>\frac{1}{3}\mu_r(D_a(r))$,
\begin{align}
\inf_{(t_0,x_0)\in Q(r)}\mP_{(t_0,x_0)}\Big(\sigma_K<\tau_{Q(\varphi_a(r))}\Big)\geq c_2\frac{\ln a}{a}.\label{H33}
\end{align}
In particular, condition  {\bf (H$_2$)} holds.
\el

\begin{proof}
Notice that
$$
\Big\{Z_{\tau_{Q(\varphi_{\sqrt{a}}(r))}}\in K\Big\}\subset\Big\{\sigma_K<\tau_{Q(\varphi_a(r))}\Big\}.
$$
It suffices to prove that there are $a_0\geq 1$ and $c_2\in(0,1)$ such that for all $a\geq a_0$ and any $(t_0,x_0)\in Q(r)$,
$$
\mP_{(t_0,x_0)}\Big(Z_{\tau_{Q(\varphi_{\sqrt{a}}(r))}}\in K\Big)\geq c_2\frac{\ln a}{a}.
$$
As $\mu_r ( \partial  Q(\varphi_{\sqrt{a}}(r)))=0$, by taking a suitable subset of $K$ if needed, we may assume without loss of generality
that $K\cap \partial  Q(\varphi_{\sqrt{a}}(r)))=\emptyset$. Then
$$
\1_K(Z_{\tau_{Q(\varphi_{\sqrt{a}}(r))}})=\sum_{0<s\leq\tau_{Q(\varphi_{\sqrt{a}}(r))}}\1_{X_{s-}\not=X_s}\1_K(Z_s).
$$
Hence by formula (\ref{e4}) and (\ref{Kappa0}), we have
\begin{align*}
\mP_{(t_0,x_0)}\Big(Z_{\tau_{Q(\varphi_{\sqrt{a}}(r))}}\in K\Big)&=\mE_{(t_0,x_0)}\int^{\tau_{Q(\varphi_{\sqrt{a}}(r))}}_0
\!\!\!\int_{\mR^d}\1_K(s,y) \kappa_s(X_s, X_s-y)\dif y\dif s\\
&\geq c_1\mE_{(t_0,x_0)}\int^{\tau_{Q(\varphi_{\sqrt{a}}(r))}}_0\!\!\!\int_{\mR^d}\frac{\1_K(s,y)}{\phi(|X_s-y|)|X_s-y|^{d}}\dif y\dif s.
\end{align*}
Since $|x-y|\leq 2|y|$ for $x\in B(\varphi_{\sqrt{a}}(r))$ and $y\notin B(\varphi_{\sqrt{a}}(r))$, by (\ref{Pr}) and definition of $\mu_r$, we have
\begin{eqnarray} \label{e:3.16}
\mP_{(t_0,x_0)}\Big(Z_{\tau_{Q(\varphi_{\sqrt{a}}(r))}}\in K\Big)
&\geq&  c_3\mE_{(t_0,x_0)}\left(\int^{\tau_{Q(\varphi_{\sqrt{a}}(r))}}_0\!\!\!
\int_{\mR^d}\frac{\1_K(s,y)}{\phi(|y|)|y|^{d}}\dif y\dif s\right) \nonumber \\
&=& c_3\mE_{(t_0,x_0)}\left(\int^{\tau_{Q(\varphi_{\sqrt{a}}(r))}}_0\!\!\!\int_{\mR^d}\frac{\1_K(s,y)}{\Phi(|y|)}\mu_r(\dif y,\dif s)\right)\\
&\geq& \frac{c_3\mu_r(K)}{\Phi(\varphi_a(r))}\mP_{(t_0,x_0)}\Big(\tau_{Q(\varphi_{\sqrt{a}}(r))}\geq\Phi(r)\Big),  \nonumber
\end{eqnarray}
where the last inequality is due to $y\in D_a(r)$ and the increasing of $\Phi$.
Lastly, by $\mu_r(K)\geq\frac{1}{3}\mu_r(D_a(r))$, (\ref{We7}) and (\ref{We6b}), we obtain that for $a\geq C^2_6\vee 4C_{8}^2=:a_0$,
\begin{align*}
\mP_{(t_0,x_0)}\Big(Z_{\tau_{Q(\varphi_{\sqrt{a}}(r))}}\in K\Big)&\geq\frac{c_3\omega_d\ln a}{6a}
\left(1-\mP_{(t_0,x_0)}\Big(\tau_{Q(\varphi_{\sqrt{a}}(r))}<\Phi(r)\Big)\right)\\
&\geq\frac{c_3\omega_d\ln a}{6a}\left(1-\frac{C_{8}}{\sqrt{a}}\right)\geq \frac{c_3\omega_d\ln a}{12a}.
\end{align*}
The proof is completed by taking $c_2=\frac{c_3\omega_d}{12}$.
\end{proof}

We can now present the

\begin{proof}[Proof of Theorem \ref{Main}]
Fix $t_0>0$ and $x_0\in\mR^d$.
\\
(i)  In this case, by assumption
$b_t (x)$ is continuous in $x$
and $b_t (x)\leq C (1+|x|)$
for all $t>0$ and $x\in \mR^d$.
Thus, by the theory of ODE, the following ODE admits at least one solution $y_t$
for $t\in [0, t_0]$:
$$
\dot y_t=-b_{t_0-t}(x_0+y_t),\ \ y_0=0.
$$
Define
$$
w(t,x):=u(t_0,x_0)-u(t_0-t,x_0+x+y_t)
$$
and
$$
\tilde b_t(x):=b_t(x+x_0+y_t)-b_t(x_0+y_t).
$$
Then
\begin{align*}
\p_t w+\sL^{\tilde b}_{t_0-t} w=0,\quad  t\in[0,t_0).
\end{align*}
Notice that by (\ref{bb0}),
$$
|\tilde b_t(x)|\leq C_b |x|/\Phi(|x|) \quad\hbox{for } |x|\leq 1.
$$
By Lemmas \ref{Le34}-\ref{Le36} and Theorem \ref{Th24}, we have
\begin{align*}
|w(t,x)|=|w(t,x)-w(0,0)|\leq 8\left(\frac{t\vee\Phi(|x|)}
{ t_0}\right)^\beta\|w\|_{L^\infty([0, t_0]\times \mR^d)}
 \quad \hbox{for } (t,x)\in Q(\Phi^{-1}(t_0)).
\end{align*}
By making the change of variables $t_0-t=t'$ and $x_0+x+y_t=x'$, and noticing that
$$
|y_t|\leq \lambda t\ \mbox{ for some } \lambda=\lambda(\|b/(1+|x|)\|_\infty)>0,
$$
we obtain the desired estimate (\ref{Es0}).
\\
\\
(ii)  In this case, define
$$
w(t,x):=u(t_0,x_0)-u(t_0-t,x_0+x).
$$
Just as above, one can conclude that (\ref{Es00}) holds.
\end{proof}

\end{document}